\documentclass[A4paper,11pt]{article}
\usepackage{latexsym}
\usepackage{mathrsfs}
\usepackage{amssymb}
\usepackage{amscd}
\usepackage[dvips]{graphicx}                  

\newtheorem{theorem}{Theorem}

\newenvironment{proof}[1][Proof]{\textbf{#1.} }{\ \rule{0.5em}{0.5em}}
\date{}
\long\def\symbolfootnote[#1]#2{\begingroup%
	\def\thefootnote{$\;$}\footnote[#1]{$^*$#2}\endgroup}
\begin{document}
	
	\title{On equivalences of polarized partition relations}
	\author{Joanna Jureczko}
\maketitle

\symbolfootnote[2]{Mathematics Subject Classification: Primary 03E02, 03E05, 03E10, 54E20.

\hspace{0.2cm}
Keywords: \textsl{partition relation, polarized  partition relation, strong sequences, product of generalized strong sequence, cardinal number, ordinal number, strongly inaccessible number}}

\begin{abstract}
	The paper deals with two notions: polarized partition relations and product of generalized strong sequences. Strong sequences were introduced by Efimov in 1965 as a usefull tool for proving famous theorems in dyadic spaces, i.e. continuous images of Cantor cube. In this paper we introduce the notion of product of generalized strong sequences and give pure combinatorial proof that  existence of product of generalized strong sequences is equivalent to polarized partition relations.
	\end{abstract}

\section{Introduction and historical background}
The study of partition relations dates back to 30' of the last century to Ramsey's paper \cite{R} however the modern notation was introduced in \cite{ER} by Erd\"os and Rado as the ordinary partition relations which concerned partitions of finite subsets of a set of a given size and the polarized partition relations which concerned partitions of finite subsets of products of sets of a given size, (where size means the cardinality of a set or order type of an ordered set. We will specify it in the concrete situations). 
The paper that led to these results was  published in 1942  by Erd\"os \cite{E} in which the author proved a generalized Ramsey theorem.

Papers that deserve attention in this topic are undoubtedly \cite{ER2, HL, W}, however, many new results have recently been proven. This shows that the topic is extremely lively and still worth exploring.

Considering the main results of our paper, at this point we will focus especially on the so-called polarized partition relations.
In the papers that we will quote below, the products of two sets are considered, e.g. two countable or two well ordered sets, of which at least one is countable. Particularly noteworthy is the result of the paper \cite{BH2} in which the authors prove that
$$\left( \begin{array}{c}
\kappa^+ \\
\kappa 
\end{array} \right) \to \left( \begin{array}{c}
\kappa \\
\kappa 
\end{array} \right)_\gamma,$$
where $\kappa$ is weakly compact and in \cite{S} for $\gamma < \kappa$ and for $\kappa$ being singular strong limit cardinal (of uncountable cofinality) which satisfies $2^\kappa > \kappa^+$ and $\gamma < cf(\kappa)$. This result belongs to Baumgartner and Hajnal, (\cite{BH}), and was known before only in the case $\gamma < \omega$.

In \cite{H}, Hajnal proved that for $\kappa$ being a measurable cardinal
$$\left( \begin{array}{c}
\kappa^+ \\
\kappa 
\end{array} \right) \to \left( \begin{array}{c}
\alpha \\
\kappa 
\end{array} \right)_{< \kappa}^{1, n},$$
where $n < \omega$ and $\alpha < \kappa^+$.

In \cite{J}, the author showed that if $\kappa^{<\kappa} = \kappa$ and there is a $\kappa$-dense ideal on $\kappa$ then
$$\left( \begin{array}{c}
\kappa^+ \\
\kappa 
\end{array} \right) \to \left( \begin{array}{c}
\alpha \\
\kappa 
\end{array} \right)^{1, 1}$$
for all $\alpha < \kappa^+$. This is true for $\kappa = \omega$, (see \cite{BH2}), or $\kappa$ being measurable, (see \cite{C}).

In \cite{GS3}, the author proved the consistency of
$$\left( \begin{array}{c}
\kappa^{++} \\
\kappa^+ 
\end{array} \right) \to \left( \begin{array}{c}
\alpha \\
\kappa^+ 
\end{array} \right)_\kappa$$
for each $\alpha \in (\kappa)^{++}$, where $\kappa$ is regular.

A lot of recent results have arisen with the assumptions of certain characteristics of the continuum. Undoubtedly, the papers \cite{GS, GS2} by Garti and Shelah deserve attention on this topic. In \cite{GS}, the authors considered the splitting number $\mathfrak{s}$ (also as a singular number) and proved that
$$\left( \begin{array}{c}
\mathfrak{s} \\
\omega 
\end{array} \right) \to \left( \begin{array}{c}
\mathfrak{s} \\
\omega 
\end{array} \right)_2.$$
In \cite{GS2}, the authors consider the consistency of ZFC with
$$\left( \begin{array}{c}
inv \\
\omega 
\end{array} \right) \to \left( \begin{array}{c}
inv \\
\omega 
\end{array} \right)_2$$
for every cardinal invariant of the continuum. Some of their results are generalized to higher cardinals. 
More considerations in these directions one can find in \cite{G} in which Garti showed independence over ZFC of polarized partitions which were posed as Problem 15 in \cite{ER2}.

The next paper which is worth attention is \cite{KW} in which the authors show other connections between  cardinal invariants and polarized partition relations. 

In most of the results in the literature we deal with the product of two sets and the partition into two colours. Few results touch upon the problem of the partition into $\omega$ colours, e.g. in \cite{L} we have
$$\left( \begin{array}{c}
\aleph_2 \\
\aleph_1 
\end{array} \right) \to \left( \begin{array}{c}
\aleph_1 \\
\aleph_1 
\end{array} \right)_\omega.$$
Moreover, the results on the polarized partition relations are mainly based on forcing methods, (see e.g. \cite{G, GS, GS2}).

In our paper we propose combinatorial methods of proving polarized relations for the special kind of system of numbers. The methods presented here rely a bit on methods presented in \cite{JJ, JJ1}. In \cite{JJ},  it is proved that the generalized Erd\"os-Rado theorem is equivalent to the existence of strong sequences. (The notion of strong sequences was introduced by Efimov in 1965). In \cite{JJ1}, the main theorem about generalized strong sequences was proved. 

Strong sequences were introducd in \cite{E3} by Efimov as a useful method for proving famous theorems in dyadic spaces (i.e. continuous images of Cantor cube). Among others, Efimov showed that there are no strong sequences in the subbase of Cantor cube although it is still an open problem in which spaces they exist. Our investigations in this topic presented in \cite{JJ3, JJ4, JJ, JJ1, JJ2} are concentrated around the possible consequences of the existence of strong sequences in spaces (without specifying what kind of space it is about). 

In the current paper, the main result concerns proving that the product of generalized strong sequences is equivalent to polarized partition relations.
We here intentionally omit polarized partition relations of the type    
$$\left( \begin{array}{c}
\alpha \\
\beta 
\end{array} \right) \to \left( \begin{array}{cc}
\gamma & \eta \\
\delta & \lambda 
\end{array} \right)$$
(for a given system of numbers $\alpha, \beta, \gamma, \delta, \eta, \lambda$)
because it is the topic of our next paper, (\cite{JJ5}), where the reader also can find more information about previous results on strong sequences.

The paper is organized as follows. In Section 2 we give basic definitions needed in further parts of the text. In Section 3 we prove the theorem on product of strong sequences. In Section 4 we prove the equivalence of Theorem 1 with polarized partition relation.

In this paper we use standard notation and terminology. For definitions and facts not cited here we refer the reader to \cite{E2, J2}.

\section{Definitions}

In the whole paper we use Greek letters to denote cardinal or ordinal numbers, (which one we will mean at the particular parts will follow from the context). 
\\\\
\noindent
\textbf{2.1.} The \textit{polarized partition relation} 
$$\left( \begin{array}{c}
\alpha_1 \\
\alpha_2 \\
...\\
\alpha_n 
\end{array} \right) \to \left( \begin{array}{c}
\beta_1 \\
\beta_2 \\
...\\
\beta_n 
\end{array} \right)_\gamma^{m_1, m_2, ..., m_n}$$
means that the following  statement is true: for all sets $A_k$ of cardinality $\alpha_k$, $(1\leqslant k \leqslant n)$, and for every function $$c \colon [A_1]^{m_1}\times [A_2]^{m_2} \times ... \times [A_n]^{m_n} \to \gamma$$ there exist $B_k \subseteq A_k$ of cardinality $ \beta_k, (1 \leqslant k \leqslant n)$, such that $$c|([B_1]^{m_1} \times [B_2]^{m_2} \times ... \times [B_n]^{m_n})$$ is constant on some $\lambda < \gamma$.
\\ 
If $m_k = 1, (1 \leqslant k \leqslant n)$ we will simply write
$$\left( \begin{array}{c}
\alpha_1 \\
\alpha_2 \\
...\\
\alpha_n 
\end{array} \right) \to \left( \begin{array}{c}
\beta_1 \\
\beta_2 \\
...\\
\beta_n 
\end{array} \right)_\gamma.$$
In the further parts of the paper we will consider infinite $\alpha_k, \beta_k$ and $\gamma$,  $(1\leqslant k \leqslant n)$.
\\\\
\noindent
\textbf{2.2.} Let $(X_k, r_k)$ be sets with two-placed relations $r_k$, $(1\leqslant k \leqslant n)$.

In the whole paper we restrict our considerations to finite products of sets, because we do not need more in this moment, but the results presented in further parts of this paper can be generalized for infinite products, (with extreme caution as is usual with infinite product operations). 
 
Let $|X_k| \geqslant \kappa_k, (1\leqslant k \leqslant n)$. Let $$X = X_1\times X_2\times ... \times X_n$$ and $\kappa = \kappa_1\cdot \kappa_2 \cdot ... \cdot \kappa_n$.

We say that $a = (a_1, a_2, ..., a_n)$ and $b = (b_1, b_2, ..., b_n) \in X$ \textit{have a bound} iff  there is $c= (c_1, c_2 ,..., c_n)$ such that $(a_k, c_k) \in r_k$ and $(b_k, c_k) \in r_k$ for every $1\leqslant k \leqslant n$.

We say that $A_k\subseteq X_k$ is \textit{$\kappa_k$-directed} if every subset of $A_k$ of cardinality less than $\kappa_k$ has a bound, $(1 \leqslant k \leqslant n)$.  

We say that $A \subseteq X$ is \textit{$\kappa$-directed} if every subset of $A$ of cardinality less than $\kappa$ has a bound.
\\\\
\noindent
\textbf{2.3.} Let $(X_k, r_k)$ be sets with relations $r_k$ and $\alpha$ and $\kappa_k, (1 \leqslant k \leqslant n)$, be cardinals. 
A sequence $(H^k_\xi)_{\xi< \alpha}$ of subsets of $X_k$ is called a \textit{$\kappa_k$-strong sequence} iff
\begin{itemize}
	\item [(1)] $H^k_\xi$ is $\kappa_k$-directed for all $\xi < \alpha$
	\item [(2)] $H^k_\xi \cup H^k_\psi$ is not $\kappa_k$-directed whenever $\xi < \psi < \alpha$, i.e. there exists $S^k_\psi \in [H^k_\psi]^{<\kappa_k}$ such that for any $\xi < \psi$ the set $H^k_\xi \cup S^k_\psi$ is not $\kappa_k$-directed. (Such a set $S^k_\psi$ is called \textit{$(k, \xi, \psi)$-destroyer}).
\end{itemize}
Let $X$ and $\kappa$ be as in Section 2.2 
and let   $H_\xi= H^1_\xi \times H^2_\xi \times ... \times H^n_\xi$.
A sequence $(H_\xi)_{\xi< \alpha}$ of subsets of $X$ is called a \textit{product $\kappa$-strong sequence} iff 
\begin{itemize}
	\item [(3)] $H_\xi$ is $\kappa$-directed for all $\xi < \alpha$,
	\item [(4)] $H_\xi \cup H_\psi$ is not $\kappa$-directed whenever $\xi < \psi < \alpha$.
	\end{itemize}
\vspace{0,5cm}
\noindent
\textbf{2.4.} Let $\Phi$ be an ordinal. A pair $(F, G)$, where $$F=(F_1, F_2, ..., F_n),$$ $$ G=(G_1, G_2, ..., G_n)$$ and $$F_k \colon X_k \to 2^\Phi$$
$$G_k \colon \Phi \to 2^{X_k}$$ is called a pair of \textit{twin functions} iff for any $k, (1 \leqslant k \leqslant n)$,  and for any $\alpha, \beta \in \Phi$ if $\alpha < \beta$ then 
 there is $a_k \in G_k (\alpha)$ such that $\beta \in F_k(a_k)$.

A map $$g=(g_1, g_2,..., g_n)$$ is called a \textit{selector of twin functions} iff for any $k, (1\leqslant k \leqslant n)$, the functions $$g_k \colon A \to X_k,\ \  A \subseteq \Phi$$ fulfill the following conditions:
\begin{itemize}
	\item [(1)] $g_k(\alpha) \in G_k(\alpha)$ for any $\alpha \in A$,
	\item [(2)] for any $\alpha, \beta \in A$, if $\alpha < \beta$ then $\beta \in F_k(g_k(\alpha))$.
\end{itemize}

\section{Theorem on product $\kappa$-strong sequences}

 Let $\beta$ and $\tau$ be cardinals. By $\beta \ll \tau$ we denote $\tau$ is \textit{strongly $\beta$-inaccessible}, i.e $\beta < \tau$ and $\alpha^\lambda < \tau$ whenever $\alpha < \tau$ and $ \lambda < \beta$. 
 
 The special case of the following theorem, (for $k=1$), was proved in \cite{JJ1}. (Theorem 1 with the very similar proof is also included in \cite{JJ5} but it is shown here intentionally because \cite{JJ5} is yet before publication).
\\

\begin{theorem}  Let $n < \omega$. For  $1\leqslant k \leqslant n$ let $\beta_k, \eta_k, \kappa_k, \mu_k$ be cardinals such that $\omega\leqslant \beta_k \ll\eta_k, \mu_k < \beta_k, \kappa_k \leqslant 2^{\mu_k}, \beta = \beta_1 \cdot \beta_2 \cdot ... \cdot \beta_n, \eta = \eta_1\cdot \eta_2 \cdot ...\cdot \eta_n, \kappa = \kappa_1\cdot \kappa_2 \cdot ...\cdot \kappa_n$  and $\beta_k, \eta_k$ be regulars. Let $X= X_1\times X_2\times ... \times X_n$ be a set and $|X_k| = \eta_k\geqslant \kappa_k$. If there exists a product $\kappa$-strong sequence $$\{H_\alpha \subseteq X \colon \alpha < \eta\}, $$ where $H_\alpha = H^1_\alpha \times H^2_\alpha \times ... \times H^n_\alpha$ and $|H^k_\alpha| \leqslant 2^{\mu_k}$ for all $\alpha < \eta$, then there exists a product $\kappa$-strong sequence $$\{T_\alpha \colon \alpha < \beta\}$$ with $|T^k_\alpha|<\kappa_k$ for all $\alpha < \beta$.
\end{theorem}

\begin{proof}
	Fix $n < \omega$.
	Let $\{H_\alpha \colon \alpha < \eta\}$  be  a product $\kappa$-strong sequence. Fix $\alpha$ and name it $\alpha_0$, (without the loss of generality one can assume that $\alpha_0 = 0$).
	For each $k, (1 \leqslant k \leqslant n)$, consider a function 
	$$f^k_{\alpha_0}\colon \eta_k\setminus \{\alpha_0\} \to [H^k_{\alpha_0}]^{< \kappa_k}$$
	such that $f^k_{\alpha_0}(\xi)=S^k_{\alpha_0}$ for some $S^k_{\alpha_0} \in [H^k_{\alpha_0}]^{< \kappa_k}$. Since $|H^k_{\alpha_0}| \leqslant 2^{\mu_k}$, hence $|[H^k_{\alpha_0}]^{< \kappa_k}|\leqslant 2^{\mu_k}< \eta_k$. It means that $f^k_{\alpha_0}$ determines a partition of $\eta_k\setminus \{\alpha_0\}$ into at most $2^{\mu_k}< \eta_k$ elements. The cardinal $\eta_k$ is regular, hence there exists $\overline{S}^k_{\alpha_0}\in [H^k_{\alpha_0}]^{< \kappa_k}$ such that $|(f^k_{\alpha_0})^{-1}(\overline{S}^k_{\alpha_0})|= \eta_k$.
	Let $$\mathcal{S}^k_{\alpha_0} = \{S^k_{\alpha_0} \in [H^k_{\alpha_0}]^{< \kappa_k} \colon |(f^k_{\alpha_0})^{-1}(S^k_{\alpha_0})|= \eta_k\}$$
	and let
	 $$\mathcal{A}^k_{\alpha_0} = \{A^k_{\alpha_0}\subseteq \eta_k\setminus \{\alpha_0\} \colon \exists(S^k_{\alpha_0}\in [H^k_{\alpha_0}]^{< \kappa_k})\  (f^k_{\alpha_0})^{-1}(S^k_{\alpha_0})= A^k_{\alpha_0} \textrm{ and } |A^k_{\alpha_0}| = \eta_k\}$$
	 be a family of pairwise disjoint sets. 
	 
	Before the continuation of the proof we make the following observation. By definition 2.3 for each $\alpha> \alpha_0$ there exists a $(k, \alpha_0, \alpha)$-destroyer. Since we consider only $k$, $(1\leqslant k \leqslant n)$,  some of them must occur at least $\beta_k$-times, $(\omega \leqslant \beta_k \ll \eta _k)$.
 Now, we are ready to continue the proof.	
	
	 For  every relevant $k$ we will construct inductively 
	 \begin{itemize}
	 	\item [a)] an increasing subsequence $\{\alpha_\gamma \colon \gamma < \beta_k\}$ of elements of  $\eta_k,$
	 	\item [b)] families $\mathcal{A}^k_{\alpha_\gamma} = \{A^k_{\alpha_\gamma}\subseteq A^k_{\alpha_{\gamma-1}}\setminus \{\alpha_\gamma\} \colon \exists(S\in [H^k_{\alpha_\gamma}]^{< \kappa_k})\  (f^k_{\alpha_\gamma})^{-1}(S)= A^k_{\alpha_\gamma} \textrm{ and } |A^k_{\alpha_\gamma}| = \eta_k\}, $
	 \end{itemize}
 where 	$$f^k_{\alpha_\gamma}\colon A \setminus \{\alpha_\gamma\}\to [H^k_{\alpha_\gamma}]^{< \kappa_k}$$
 such that $A \in \mathcal{A}^k_{\alpha_{\gamma-1}}$ and $f^k_{\alpha_\gamma}(\xi) = S$ for some $(k, \alpha_\sigma, \alpha_\gamma)$-destroyers, whenever $\alpha_\sigma < \alpha_\gamma$.
 
 Assume that we have constructed increasing subsequence $\{\alpha_\gamma \colon \gamma < \beta_k\}$ of $\eta_k$ and families $\mathcal {A}^k_{\alpha_\gamma}$ as was done above. 
 
 Next, choose $\alpha > \alpha_\gamma, (\alpha < \eta_k)$ such that there exists a $(k, \alpha_\gamma, \alpha)$-destroyer $S^k_{\alpha} \in [H^k_{\alpha}]^{< \kappa_k}$ and denote this $\alpha$ by $\alpha_\delta$, where $\delta= \gamma+1$.
 For each $A^k_{\alpha_\gamma} \in \mathcal{A}^k_{\alpha_\gamma}$ define a function 
 $$f^k_{\alpha_\delta}\colon A^k_{\alpha_\gamma} \setminus \{\alpha_\delta\}\to [H^k_{\alpha_\delta}]^{< \kappa_k}$$
 such that $f^k_{\alpha_\delta}(\xi) = S^k_{\alpha_\delta}$ for some $(k, \alpha_\gamma, \alpha_\delta)$-destroyer $S^k_{\alpha_\delta} \in [H^k_{\alpha_\delta}]^{\kappa_k}$ (and any $\alpha_\gamma < \alpha_\delta$).
 Since $|H^k_{\alpha_\delta}| \leqslant 2^{\mu_k}< \eta_k$ and $|A^k_{\alpha_\gamma}| = \eta_k$ the function $f^k_{\alpha_\delta}$ determines a partition of $A^k_{\alpha_\gamma} \setminus \{\alpha_\delta\}$ into at most $2^{\mu_k}$ elements. Hence there exists a $(k, \alpha_\gamma, \alpha_\delta)$-destroyer $\overline{S}^k_{\alpha_\delta} \in [H^k_{\alpha_\delta}]^{< \kappa_k}$ such that $|(f^k_{\alpha_\delta})^{-1}(\overline{S}^k_{\alpha_\delta})|= \eta_k$. Let 
 $$\mathcal{S}^k_{\alpha_\delta} = \{S^k_{\alpha_\delta} \in [H^k_{\alpha_\delta}]^{< \kappa_k} \colon |(f^k_{\alpha_\delta})^{-1}(S^k_{\alpha_\delta})|= \eta_k\}$$
 and let
 $$\mathcal{A}^k_{\alpha_\delta} = \{A^k_{\alpha_\delta}\subseteq  A^k_{\alpha_\gamma} \setminus \{\alpha_\delta\}\colon \exists(S^k_{\alpha_\delta}\in [H^k_{\alpha_\delta}]^{< \kappa_k})\  (f^k_{\alpha_\delta})^{-1}(S^k_{\alpha_\delta}) =A^k_{\alpha_\delta} \textrm{ and }|A^k_{\alpha_\delta}|=\eta_k\}$$
 be a family of pairwise disjoint sets.
 
 If $\delta$ is limit, we consider 
 $$f^k_{\alpha_\delta}\colon \bigcap_{\rho< \gamma }A^k_{\alpha_\rho} \setminus \{\alpha_\delta\}\to [H^k_{\alpha_\delta}]^{< \kappa_k}$$
 	for $\bigcap_{\rho< \gamma }A^k_{\alpha_\rho} \not = \emptyset$.
 	(Note that we can always find such a non-empty intersection because at each step of the induction we divide each subset of $A^k_{\alpha_\delta} \in \mathcal{A}^k_{\alpha_\delta}$ for $\delta < \beta_k$ into pairwise disjoint sets and use elements (of cardinality $\eta_k$) of such obtained partition in the next step of induction).
 		The induction step is complete.
 		
 		Now, define a sequence 
 		$$T_{\alpha_\gamma} = T^1_{\alpha_\gamma} \times T^2_{\alpha_\gamma} \times ... \times T^n_{\alpha_\gamma}$$ for all $\gamma< \beta$ in the following way:
 \[ T^k_{\alpha_\gamma}=
 \left\{\begin{array}{rl}
 S^k_{\alpha_\gamma} \in  [H^k_{\alpha_\gamma}]^{< \kappa_k} & \mbox{if }\exists k\  S^k_{\alpha_\gamma} \in \mathcal{S}^k_{\alpha_\gamma}\\
P^k_{\alpha_\gamma} \in  [H^k_{\alpha_\gamma}]^{< \kappa_k} & \mbox{otherwise},
 \end{array}\right. \]
 where $P^k_{\alpha_\gamma}$ is arbitrarily chosen. Thus, we have defined at least one product $\kappa$-strong sequence $\{T_{\alpha_\gamma} \colon \alpha_ \gamma < \beta\}$ of the required property.
 
 Suppose now, that at least one of product $\kappa$-strong sequences $$\{T_{\alpha_\gamma} \colon \alpha_ \gamma < \beta\}$$ has length $\zeta> \beta$, i. e. there exists $k, (1\leqslant k \leqslant n),$ which occurs $\zeta$ times, i. e. there is a sequence $$\{S^k_{\alpha_\gamma} \colon \alpha_ \gamma < \zeta\}$$
 such that $S^k_{\alpha_\gamma} \in \mathcal{S}^k_{\alpha_\gamma}$.
 By our construction, each $S^k_{\alpha_\gamma}$ determines a set $A \in \mathcal{A}^k_{\alpha_\gamma}$ such that $|A| = \eta_k$.
 Let $$\nu_k = \sup\{|A|\colon  A \in \mathcal{A}^k_{\alpha_\gamma}, \gamma < \zeta\}.$$
 Then, there would exist $\nu_{k} > \eta_k$ pairwise disjoint sets $A \in \mathcal{A}^{k}_{\alpha_\gamma},$ where $ |A| = \eta_{k}$. A contradiction.
\end{proof}

\section{Main results}
In this section we prove that Theorem 1 is equivalent to the following theorem.

\begin{theorem}
Let $n < \omega$. For  $1\leqslant k \leqslant n$ let $\beta_k,\eta_k, \kappa_k,  \mu_k$ be cardinals such that $\omega\leqslant \beta_k \ll\eta_k, \mu_k < \beta_k, \kappa_k \leqslant 2^{\mu_k},\beta = \beta_1 \cdot \beta_2 \cdot ... \cdot \beta_n, \mu = \mu_1\cdot \mu_2 \cdot ...\cdot \mu_n$  and $\beta_k, \eta_k$ be regulars. Then
	$$\left( \begin{array}{c}
	\eta_1 \\
	\eta_2 \\
	...\\
	\eta_n 
	\end{array} \right) \to \left( \begin{array}{c}
	\beta_1 \\
	\beta_2 \\
	...\\
	\beta_n 
	\end{array} \right)_\mu,$$
	i. e. for each colouring function 
	$c \colon X \to \mu$
	such that $X= X_1\times X_2\times ... \times X_k$ with $|X_k| = \eta_k \geqslant \kappa_k$ there exists $A = A_1\times A_2\times ... \times A_k $ such that $A_k\subseteq X_k$ with $|A_k| = \beta_k$ and $c^{-1}(\lambda) = A$ for some  $\lambda < \mu$, ($c$ is constant on some set $A$ of cardinality $\beta$). 
\end{theorem}

For this purpose we use Theorem 3 (below) as an auxiliary result. Thus the proving scheme will be as follows. We firstly show that Theorem 1 is equivalent to Theorem 3 and secondly we show that Theorem 2 is equivalent to Theorem 3.

\begin{theorem} Let $n < \omega$. For  $1\leqslant k \leqslant n$ let $\beta_k, \eta_k, \kappa_k, \mu_k$ be cardinals such that $\omega\leqslant \beta_k \ll\eta_k, \mu_k < \beta_k, \kappa_k \leqslant 2^{\mu_k}, \beta = \beta_1 \cdot \beta_2 \cdot ... \cdot \beta_n,  \eta = \eta_1\cdot \eta_2 \cdot ...\cdot \eta_n$  and $\beta_k, \eta_k$ be regulars. Let $X= X_1\times X_2\times ... \times X_n$ be a set and $|X_k| = \eta_k\geqslant \kappa_k$. If there exists a pair $(F, G)$ of twin functions, where
	$$F=(F_1, F_2, ..., F_n),\ \ \  G=(G_1, G_2, ..., G_n)$$ and $$F_k \colon X_k \to 2^\eta,\ \ \ G_k \colon \eta \to 2^{X_k}$$
	  such that $|G_k(\alpha)|\leqslant 2^{\mu_k}$ for every $\alpha < \eta_k$, $(1\leqslant k \leqslant n)$, then there exists a selector  of twin functions $g=(g_1, g_2, ..., g_n)$ such that $g_k\colon A \to X_k$ with $A \subseteq \eta$ and $|A| = \beta$. 
	\end{theorem}

\begin{proof}  Fix $n < \omega$.
	\newline

(\textbf{Theorem 1 implies Theorem 3})
Take a pair $(F, G)$ of twin functions such that $|G_k(\alpha)|\leqslant 2^{\mu_k}$ for any $k, (1\leqslant k \leqslant n)$.
For any $\alpha < \eta_k$ consider
$$H^k_\alpha = G_k(\alpha) \setminus (F_k)^{-1}(\alpha).$$
	Notice that $$\{H_\alpha \colon \alpha < \eta\}$$ is a product $\kappa$-strong sequence. If not, then by definition 2.3 the set $H^k_\alpha\cup H^k_\gamma$ would be $\kappa_k$-directed for any $\alpha< \gamma$ and any $k, (1\leqslant k \leqslant n)$. By definition 2.4 for all $k$ there exists  $a_k \in H^k_\gamma$ such that $\gamma \in F_k(a_k)$. Hence  $a_k \in (F_k)^{-1}(\gamma)$ for any $k, (1\leqslant k \leqslant n)$ which contradicts to definition of $H^k_\gamma$.
		
		By Theorem 1, there  exists a product $\kappa$-strong sequence $$\{T_\alpha \colon \alpha < \beta\}$$ such that $T^k_\alpha \in [H^k_\alpha]^{<\kappa_k}$, (see the proof of Theorem 1),  for any $\alpha < \beta$. For any $k, (1\leqslant k \leqslant n)$, consider a function $$g_k\colon A \to X_k$$
		such that $A \subseteq \eta, |A|=\beta$ and $g_k(\alpha) \in T^k_\alpha$ for any $k$, $(1\leqslant k \leqslant n),$ and $\alpha \in A$. Obviously, $g_k(\alpha) \in G_k(\alpha)$, hence condition $(1)$ in definition 2.4 is fulfilled. Moreover, by definition 2.4, since $(F, G)$ are twin function we have that $$\gamma \in  F_k(g_k(\alpha))$$ for any $\alpha < \gamma$ and any $k, (1\leqslant k \leqslant n)$. Hence condition $(2)$ of  definition 2.4 is fulfilled. 
\\

(\textbf{Theorem 3 implies Theorem 1})
Assume that there exists a product $\kappa$-strong sequence $$\{H_\alpha \colon \alpha < \eta\}$$ such that $|H^k_\alpha|\leqslant 2^{\mu_k}$ for any $k, (1\leqslant k \leqslant n)$. For any $k, (1\leqslant k \leqslant n)$, consider the set
$$\mathcal{C}^k_\alpha = \{T \in [H^k_\alpha]^{<\kappa_k} \colon  T \cup H^k_\xi \textrm{ is not $\kappa_k$-directed for any $\xi < \alpha$}\}.$$
By definition 2.3,  $\mathcal{C}^k_{\alpha} \not = \emptyset$ for some $k$ and $\alpha < \eta$.
For such $k$ consider 
$$\mathcal{D}_k = \{T \colon T \in \mathcal{C}^k_{\alpha} \textrm{ for some } \alpha < \eta\}$$ and define a pair of  functions $F, G$, where
$F=(F_1, F_2, ..., F_n),$ $G=(G_1, G_2, ..., G_n)$ and
$$F_k \colon \mathcal{D}_k \to 2^\eta \textrm{ such that }F_k(T) =\{\xi < \eta \colon T \cup H^k_\xi \textrm{ is not $\kappa_k$-directed}\}$$
$$G_k \colon \eta \to 2^{\mathcal{D}_k} \textrm{ such that }G_k(\alpha)= \mathcal{C}^k_\alpha.$$  We will show that $(F,G)$ are twin functions. 

Let $\alpha < \gamma$. Then for some $T \in G_k(\alpha)$ the set $T \cup H^k_\gamma$ is not $\kappa_k$-directed. Hence $\gamma \in F_k(T)$. The deifinition 2.4 is fulfilled. By Theorem 3, there exists a selector of twin functions $g = (g_1, g_2,..., g_n)$ such that $$g_k \colon A \to \mathcal{D}_k$$ where $A \subseteq \eta$,$|A| = \beta$ an $k, (1\leqslant k \leqslant n)$.

 By $(1)$ in definition 2.4 and the above construction, $g_k(\alpha) \in G_k(\alpha)$ for any $\alpha \in A$, i.e. $g_k(\alpha) = T$ for some $T \in \mathcal{C}^k_\alpha$. Hence $|g_k(\alpha)|< \kappa_k$. Moreover, since $T \in [H^k_\alpha]^{< \kappa_k}$ the set $g_k(\alpha)$ is $\kappa_k$-directed.
 
 Now by $(2)$ in definition 2.4, for $\gamma \in A$ we have $\gamma \in F_k(g_k(\alpha))$. By definition of $F_k$ the set $g_k(\alpha) \cup H^k_\gamma$ is not $\kappa_k$-directed for any $\gamma \in F_k(g_k(\alpha))$. 
 
 Since each $g_k$ determines some $T \in [H^k_\alpha]^{<\kappa_k}$ of required properties hence 
 $\{g(\alpha) \colon \alpha \in A\}$ is the required product $\kappa$-strong sequence.
 \\
 
(\textbf{Theorem 2 implies Theorem 3})
Enumerate $X_k = \{x^k_{\alpha} \colon \alpha < \eta_k\}$, $(1 \leqslant k \leqslant n)$. Let $F, G$ be functions such that $$F=(F_1, F_2, ..., F_n),\ \ \  G=(G_1, G_2, ..., G_n)$$
and $$F_k \colon X_k \to 2^{\eta}\ \ \ G_k\colon \eta \to 2^{X_k}.$$
Without loss of generality we can assume that for any $\xi < \eta$ the the set $G_k(\xi)$ is of the form $$G_k(\xi) = \{x^k_\delta (\xi) \colon \delta < \lambda_k\}$$
for some $\lambda_k \leqslant 2^{\mu_k}$.
Let $X =X_1\times X_2\times ... \times X_n$ and $$c \colon X  \to \mu$$  be a colouring function such that if $c(x^1_\gamma(\alpha_1), x^2_\gamma(\alpha_2), ..., x^n_\gamma(\alpha_n)) = \gamma$ then $\alpha_m \in F_k(x^k_{\gamma}(\alpha_k))$ for some $\alpha_m> \alpha_k$. It is easy to check that $(F, G)$ are twin functions.  

By Theorem 2, there exists a set $B_1\times B_2\times ...  \times B_n = c^{-1}(\gamma)$ for some $\gamma <\mu$ and $B_k \subseteq X_k$ such that $|B_k| = \beta_k$, $(1 \leqslant k \leqslant n)$.

Consider a function $g_k\colon A \to X_k$ such that $A \subseteq\eta$ and  $g_k(\alpha_k) = x^k_\gamma(\alpha_k)$ for any $\alpha_k \in A$.
By the above construction, $g$ is a selector of the twin functions $(F, G)$. 
\\

(\textbf{Theorem 3 implies Theorem 2})
Let $X =X_1\times X_2\times ... \times X_n$ and $$c \colon X \to \mu$$ be a colouring function.
For each $\gamma < \mu$ take $c^{-1}(\gamma)$. Then
 $$X = \{A_\gamma \colon A_\gamma = c^{-1}(\gamma), \gamma < \mu\}$$
 is a partition of $X$.
 We will show that there exists $\gamma_0 < \mu$ such that $|A_{\gamma_0}| =\beta$, (i.e. $c$ is constant on some set $ A_{\gamma_0}$, $|A_{\gamma_0}| = \beta$, $\gamma_0 < \mu$).
 
 In order to do this we will define functions
$$F=(F_1, F_2, ..., F_n),\ G=(G_1, G_2, ..., G_n)$$ such that for any $k$, $(1\leqslant k \leqslant n)$,
$$F_k \colon X_k \to 2^\eta\ \ \ G_k \colon \eta \to 2^{X_k}$$ 
and 
$$F_k(x^k_\gamma(\alpha_k)) = \{\alpha_m\colon  \alpha_m> \alpha_k, (x^1_\gamma(\alpha_1), x^2_\gamma(\alpha_2), ..., x^n_\gamma(\alpha_n))\in A_\gamma, \ \alpha_1< \alpha_2< ...<\alpha_n\}$$
and 
$$G_k(\alpha_k) = \{x^k_\gamma(\alpha_k)\colon (x^1_\gamma(\alpha_1), x^2_\gamma(\alpha_2), ..., x^n_\gamma(\alpha_n))\in A_\gamma, \  \alpha_1< \alpha_2< ...<\alpha_n\}.$$  
We will show that $(F, G)$ are twin functions.
To show this take $x^k_\gamma (\alpha_k) \in G_k(\alpha_k)$. Then $(x^1_\gamma(\alpha_1), x^2_\gamma(\alpha_2), ..., x^n_\gamma(\alpha_n))\in A_\gamma$. Then $\alpha_m \in F_k(x^k_\gamma(\alpha_k))$ for any $\alpha_m> \alpha_k$. 

By Theorem 3, there exists a selector of twin functions $g=(g_1, g_2, ..., g_n)$, where
$$g_k \colon B \to X_k,$$ 
$B \subseteq \eta, |B| = \beta$ and $g_k(\alpha_k) = x^k_\gamma(\alpha_k)$ for some $\gamma<\mu$. It measn that for all $\alpha_1, \alpha_2, ..., \alpha_n \in B$ if  $\alpha_1< \alpha_2< ...<\alpha_n$ then $(x^1_\gamma(\alpha_1), x^2_\gamma(\alpha_2), ..., x^n_\gamma(\alpha_n))\in A_\gamma$. 

Suppose that $g(B) \not = A_\gamma$ for any $\gamma < \mu$.
Take
$$(x^1_\gamma(\alpha_1), x^2_\gamma(\alpha_2), ..., x^n_\gamma(\alpha_n))\in A_\gamma,$$
where $\alpha_1< \alpha_2< ...<\alpha_n$ and  
$$(x^1_\delta(\beta_1), x^2_\delta(\beta_2), ..., x^n_\delta(\beta_n))\in A_\delta,$$
where $\beta_1< \beta_2< ...<\beta_n$ and $\gamma \not = \delta$.
By definition 2.4, 
$g_k(\alpha_k) = x^k_\gamma(\alpha_k)$ and hence $\alpha_m \in F_k(g_k(\alpha))$ for $\alpha_m> \alpha_k$
and 
$g_k(\beta_k) = x^k_\delta(\beta_k)$ and hence $\beta_m \in F_k(g_k(\beta_k))$ for $\beta_m> \beta_k$.
Now, if $\alpha_m< \beta_m$ then
$\alpha_m \in F_k(g_k(\alpha_k))$ and then 
$$(g_1(\alpha_1), g_2(\alpha_2), ..., g_n(\alpha_n)) \in A_\gamma$$
and
$$(g_1(\beta_1), g_2(\beta_2), ..., g_n(\beta_n)) \in A_\gamma.$$
But, if  $\beta_m< \alpha_m$ then
$\beta_m \in F_k(g_k(\beta_k))$ and then 
$$(g_1(\alpha_1), g_2(\alpha_2), ..., g_n(\alpha_n)) \in A_\delta$$
and
$$(g_1(\beta_1), g_2(\beta_2), ..., g_n(\beta_n)) \in A_\delta.$$
Since $A_\gamma \cap A_\delta = \emptyset$, whenever $\gamma \not = \delta$, we  obtain a contradiction.

Hence $g(B) = A_{\gamma_0}$ for some $\gamma_0 < \mu$. 
To complete the proof it is enough to show that $g_k(\alpha_k) \not = g_k(\beta_k)$ for any $k$, $(1\leqslant k\leqslant n)$, whenever $\alpha_k\not = \beta_k$. 
\\Indeed. Consider the case $\alpha_k< \alpha_m < \beta_k<\beta_m$. If $g_k(\alpha_k) = g_k(\beta_k)$ then 
$$\alpha_m \in F_k(g_k(\alpha_k)) = F_k(g_k(\beta_k))$$
and 
$$\beta_m \in F_k(g_k(\alpha_k)) = F_k(g_k(\beta_k))$$ which contradicts with definition 2.4. 
Thus $|A_{\gamma_0}| = \beta$.
\end{proof}

\section{Results for singulars}

 Shelah, in \cite{S}, showed that
$$\left( \begin{array}{c}
\kappa^+ \\
\kappa 
\end{array} \right) \to \left( \begin{array}{c}
\kappa +1 \\
\kappa 
\end{array} \right)_\gamma,$$
also holds for a singular strong limit cardinal $\kappa$ which satisfies $2^\kappa > \kappa^+$ and $\gamma < cf(\kappa)$.

Also Theorem 1 and Theorem 3 can be proved for singular numbers with using the appropriately modified proofs presented in Section 3 and Section 4. We will therefore conclude the paper with formulating theorems analogous to Theorem 1 Theorem 3 in the singular case.

\begin{theorem} Let k be a natural number such that $1\leqslant k \leqslant 2$.
	Let $\beta_k, \eta_k, \kappa_k, \mu_k$, be cardinals such that
	\begin{itemize}
		\item $\omega \leqslant \beta_k \ll\eta_k$, $\mu_k < \beta_k, \kappa_k \leqslant 2^{\mu_k},$
		\item $\beta_1, \beta_2, \eta_2$ are singular 
		\item $\eta_1 = \eta_2^+$ , $\beta_2 = \eta_2$, $\beta_1 = \beta_2+1$
		\item $2^{\beta_k} >\beta^+_k$, $\mu_k < cf(\beta_k)$.
	\end{itemize}
	Let $X = X_1 \times X_2$ be a set such that $|X_k| = \eta_k\geqslant \kappa_k$ and $\kappa = \kappa_1 \cdot \kappa_2, \eta = \eta_1 \cdot \eta_2 , \beta = \beta_1 \cdot \beta_2$. If there exists a product $\kappa$-strong sequence $$\{H_\alpha\subseteq X \colon \alpha < \eta\}$$ such that $|H^k_\alpha| \leqslant 2^{\mu_k}$ for all $\alpha < \eta$ then there exists a product $\kappa$-strong sequence $$\{T_\alpha  \colon \alpha < \beta\}$$ with $|T^k_\alpha| < \kappa_k$ for all $\alpha < \beta$.
	\end{theorem}

\begin{theorem}
	Let k be a natural number such that $1\leqslant k \leqslant 2$.
	Let $\beta_k, \eta_k, \kappa_k, \mu_k$, be cardinals $\beta = \beta_1\cdot \beta_2$, $\eta = \eta_1\cdot \eta_2$ and
	\begin{itemize}
		\item $\omega \leqslant \beta_k \ll\eta_k$, $\mu_k < \beta_k, \kappa_k \leqslant 2^{\mu_k},$
		\item $\beta_1, \beta_2, \eta_2$ are singular 
		\item $\eta_1 = \eta_2^+$ , $\beta_2 = \eta_2$, $\beta_1 = \beta_2+1$
		\item $2^{\beta_k} >\beta^+_k$, $\mu_k < cf(\beta_k)$.
	\end{itemize}
		If there exists a pair $(F, G)$ of twin functions, where 
$$F=(F_1, F_2, ..., F_n),\ \ \  G=(G_1, G_2, ..., G_n)$$
and $$F_k \colon X_k \to 2^{\eta}\ \ \ G_k\colon \eta \to 2^{X_k}$$
and $|X_k|=\eta_k \geqslant \kappa_k$				
		such that $|G_k(\alpha)|\leqslant 2^{\mu_k}$ for every $\alpha < \eta$ then there exists a selector of twin functions $g = (g_1, g_2)$, $(g_k \colon A \to X_k)$,   such that  $A \subseteq \eta$ and $|A| = \beta$. 
\end{theorem}	 

\noindent
\textbf{Acknowledgments} The author would like to thank the reviewer for incredible patience,  detailed and insightful reading of the text and for all the comments   that undoubtedly helped to improve the text and avoid inaccuracies and omissions.

\noindent
{\sc Joanna Jureczko}
\\
Wroc\l{}aw University of Science and Technology, Wroc\l{}aw, Poland
\\
{\sl e-mail: joanna.jureczko@pwr.edu.pl}

\end{document}